\documentclass{article}
\usepackage{amssymb}

\usepackage{graphicx}
\usepackage{amsmath}
\usepackage{a4wide}

\input{tcilatex}

\begin{document}

\title{Equivariant Cohomology and Representations of the Symmetric Group}
\author{Michael Atiyah}
\date{}
\maketitle

\thispagestyle{empty}

\section{Introduction}

In a recent paper \cite{1} I\ have shown how to construct a continuous map 
\begin{equation*}
f:C_{n}(R^{3})\rightarrow U(n)/T^{n}
\end{equation*}
from the configuration space of $n$ ordered distinct points of $R^{3}$ to
the flag manifold of $U(n)$ which is compatible with the natural action of
the symmetric group $\Sigma _{n}$ on both spaces. \ I also noted in \cite{1}
that the action of $\Sigma _{n}$ on the rational cohomology of either space
coincides with the regular representation but that the homomorphism $f^{\ast
}$ induced by $f$ cannot possibly be an isomorphism. \ In fact the
Poincar\'{e} polynomials of the two spaces are quite different. 
\begin{eqnarray}
P_{t}(C_{n}(R^{3})) &=&(1+t^{2})(1+2t^{2})...(1+(n-1)t^{2})  \notag \\
&&  \TCItag{1.1} \\
P_{t}(U(n)/T^{n}) &=&(1+t^{2})(1+t^{2}+t^{4})...(1+t^{2}+t^{4}+...+t^{2n-2})
\notag
\end{eqnarray}
although, for $t=1,$ both yield $n!.$ \ This shows that the various
components of the regular representation acquire quite different gradings in
the two cases and $f^{\ast }$ can only pair off some of these components.

The significance of the map $f$ for the comparison of the cohomology of the
two spaces appears therefore to be disappointingly limited. \ However the
map $f$ of \cite{1} has a further naturality property which we have so far
not exploited, namely that it is also compatible with the natural action of
the rotation group $SO(3)$ on $P(C^{n}).$ \ We can therefore consider the $%
SO(3)$-equivariant cohomology of both spaces and the homomorphism induced by 
$f^{\ast }$ in this context. \ As we shall see the situation is now totally
different and $f^{\ast }$ becomes much more interesting.

Let us first recall that for any compact Lie group $G$ and any $G$-space $X$
we define the $G$-equivariant cohomology $H_{G}^{\ast }(X)$ to be the
ordinary cohomology of the space $X_{G}$ which is the fibering over the
classifying space $BG$ with fibre $X.$ \ It follows that $H_{G}^{\ast }(X)$
is a module over the ring 
\begin{equation*}
H_{G}^{\ast }\text{(point)}=H^{\ast }(BG)
\end{equation*}
(we can take any ring of coefficients, for example the rationals $Q).$ \ For 
$G=SO(3)$ and coefficients $Q,$ $H_{G}^{\ast }$(point) is a polynomial ring
in one generator $u$ of degree $4$ (the Pontrjagin class).

The main result of this paper is then\medskip

\noindent \textbf{Theorem 1.\quad The homomorphism} 
\begin{equation*}
f^{\ast }:H_{SO(3)}^{\ast }(U(n)/T^{n},Q)\rightarrow H_{SO(3)}^{\ast
}(C_{n}(R^{3}),Q)
\end{equation*}

\noindent \textbf{is injective and its cokernel is annihilated by a power of 
}$u.\medskip $

This theorem can be put into more concrete form. \ First we shall see that
both sides are free modules over $Q[u]$ of $\func{rank}n!,$ and that they
have canonical bases. \ This identifies them both with the regular
representation of $\Sigma _{n}$ over $Q[u].$ \ Thus $f^{\ast }$ can be
explicitly represented by a matrix $A(u)$ of polynomials. \ Theorem 1
asserts that $\det A(u)\neq 0.$

For $u=0,$ we will have $\det A(0)=0,$ so that $A(0)$ is a singular matrix.
\ In fact $A(0)$ represents the induced homomorphism $f^{\ast }$ in ordinary
cohomology, which pairs off cohomology in both spaces in the same dimension.
\ On the other hand the coefficient of some other power $u^{k}$ will give a
homomorphism between the ordinary cohomology group which \textbf{lowers
dimension }by $4k.$ \ Thus, by using all these dimension-lowering operators, 
$f^{\ast }$ effectively matches up all the cohomology on both sides (as $%
\Sigma _{n}$-modules).

Consider, for example, the case $n=3,$ so that our Poincar\'{e} polynomials
are

\begin{equation*}
\begin{array}{lll}
P_{t}(C_{3}(R^{3})) & = & (1+t^{2})(1+2t^{2})=1+3t^{2}+2t^{4} \\ 
&  &  \\ 
P_{t}(U(3)/T^{3}) & = & (1+t^{2})(1+t^{2}+t^{4})=1+2t^{2}+2t^{4}+t^{6}
\end{array}
\end{equation*}
To match these up we see that we must have an operation that lowers degree
by $4$ which matches up $t^{6}$ in the flag manifold with one of the $t^{2}$
terms in the configuration space. \ This means that our matrix $A(u)$ will
look like 
\begin{equation*}
A(u)=A_{0}+A_{1}u
\end{equation*}
with $A_{0}$ being the singular matrix (of $\func{rank}5)$ describing
ordinary cohomology and $A_{1}$ being the shift operator disposing of the
extra term in dimension $6.$

More generally, for any $n,$ we can compare the two Poincar\'{e} series
given in (1.1). \ Since they both involve a factor $(1+t^{2})$ and they
agree for $t^{2}=1$ it follows that their difference 
\begin{equation*}
\phi
(t)=(1+t^{2})(1+2t^{2})...(1+(n-1)t^{2})-(1+t^{2})(1+t^{2}+t^{4})...(1+t^{2}+...+t^{2(n-1)})
\end{equation*}
is divisible by $(1+t^{2})(1-t^{2})=(1-t^{4}),$ so that $\phi
(t)=(1-t^{4})\psi (t)$ with $\psi (t)$ a polynomial. \ This is consistent
with Theorem 1. \ In fact, if we use the natural total grading of
equivariant cohomology, then $\psi (t)=\dfrac{\phi (t)}{1-t^{4}}$ is the
difference of the Poincar\'{e} series of the two free modules and by Theorem
1 this must be the Poincar\'{e} polynomial of the cokernel of $f^{\ast }.$ \
Notice that $\psi (t)$ has non-negative coefficients.

There are a number of variants of Theorem 1 which may be useful. \ We can
replace $SO(3)$ by its maximal torus $SO(2)$ or we can lift to Spin$(3)$. \
We can also work with integer cohomology, although we may have problems with
torsion at certain stages. \ More interestingly we can replace cohomology by 
$K$-theory. \ In this case our coefficient ring for equivariant $K$-theory
will be a ring of Laurent polynomials. \ We shall examine this case
carefully in \S 4.

The map $f$ of $\cite{1}$, or a variant of it, can be generalised so that $%
U(n)/T^{n}$ gets replaced by $G/T$ for any compact simple Lie group and $%
\Sigma _{n}$ gets replaced by the Weyl group. \ This follows from results of
Bielawski \cite{4} and will be explained elsewhere. \ We can then extend the
results of this paper to this more general case.

The interest of the whole story lies in the dimension shifting operators
that enter in matching up representations of the Weyl group. \ This is
reminiscent of what happens in the construction of the Springer
representation using the theory of perverse sheaves. \ The use of
equivariant $K$-theory also resembles some of the work of Lusztig \cite{6}
in this area. \ These ideas might merit further exploration.

I should also point out that, although I\ have focused on the flag manifold $%
G/T,$ there are parallel results for the other homogeneous spaces that occur
as adjoint orbits of $G.$ \ I shall discuss this (for $G=U(n))$ in \S 3.

The proof of Theorem 1 is in fact quite routine and will be explained in the
next section. \ Once the idea has occurred that equivariant cohomology might
be interesting the proof presents no difficulties.

\section{The fixed-point formula}

We shall first review some basic facts concerning equivariant cohomology
(see \cite{2} or \cite{5} for a detailed exposition) and for simplicity we
shall concentrate on the case where the group is the circle, which we shall
denote by $S.$ \ Then, for any $G$-space $X,\;H_{S}^{\ast }(X)$ is a module
over 
\begin{equation*}
\Lambda =H_{S}^{\ast }(pt)=H^{\ast }(BS)=H^{\ast }(CP_{\infty })
\end{equation*}
which is a polynomial ring with one generator $t$ of dimension $2$ (the
ground ring can be $Z$ or $Q$, but we shall focus on $Q).$

Let $F\subset X$ be the subspace of fixed points of $S,$ then the
``localisation theorem'' asserts that 
\begin{equation}
H_{S}^{\ast }(X)\rightarrow H_{S}^{\ast }(F)\text{\textbf{is an isomorphism
modulo modules annihilated by some power of} }t.  \tag{2.1}
\end{equation}
Here some mild assumptions need to be imposed on $X,$ which are certainly
satisfied in our cases when $X$ is a compact manifold (or the complement of
a finite member of closed submanifolds in a compact manifold).

The proof of (2.1) follows from a general spectral sequence argument, but it
can also be done more simply by a Meyer-Vietoris argument. \ The key point
is that, if $F$ is empty, so that all orbits of $S$ in $X$ are $1$%
-dimensional, then \textbf{over }$Q,\;H_{S}^{\ast }(X)\cong H^{\ast }(X/S)$
is finite-dimensional and hence annihilated by $t^{m}$ provided $2m>\dim
(X/S).$

Now since $S$ acts trivially on $F$ we have 
\begin{equation}
H_{S}^{\ast }(F)\cong H^{\ast }(F)\otimes _{Q}\Lambda  \tag{2.2.}
\end{equation}
so that $H_{S}^{\ast }(F)$ is the free $\Lambda $-module generated by a
basis of $H^{\ast }(F).$

Suppose in addition that $X$ has \textbf{no odd dimensional cohomology},
then the spectral sequence of the fibration $X_{S}\rightarrow BS$ shows that 
$H_{S}^{\ast }(X)$ is the free $\Lambda $-module generated by any lift of a
basis of $H^{\ast }(X)$ to $H_{S}^{\ast }(X).$

We shall now apply these general facts in our special situation. \ First we
take $X=C_{n}(R^{3})$ and we take $S=SO(2)\subset SO(3)$ to be given by
rotation about a preferred axis in $R^{3}.$ \ Since $X$ has no
odd-dimensional cohomology \cite{1} it follows that $H_{S}^{\ast }(X)$ is a
free $\Lambda $-module. \ We shall return later to consider how to construct
an explicit basis by lifting generators of $H^{\ast }(X).$ \ The fixed
points of the action of $S$ are just given by 
\begin{equation*}
F=C_{n}(R^{1})\subset C_{n}(R^{3})
\end{equation*}
where $R^{1}$ is the axis of the rotation. \ But a point in $F$ is just a
permutation of an increasing sequence 
\begin{equation*}
x_{1}<x_{2}<...<x_{n}
\end{equation*}
so that $F$ consists of $n!$ contractible components indexed by $\Sigma
_{n}. $ \ In fact $F$ can be identified with the regular points in the Lie
algebra of the maximal torus of $U(n),$ the components being the Weyl
chambers (this aspect generalises to other Lie groups $G$ instead of $U(n)).$
\ Hence $H_{S}^{\ast }(F)$ is a free $\Lambda $-module of rank $n!.$ \ 
\textbf{It can be viewed as the regular representation of }$\Sigma _{n}$ 
\textbf{over }$\Lambda .$

Next we take $X=U(n)/T^{n}.$ \ Since $S$ acts irreducibly on the projective
space $P(C^{n})$ it has just $n$ fixed points in this space (the
eigenvectors) and hence the only fixed flags are those got by permuting the
eigenvectors. \ Hence there are $n!$ fixed points, indexed by $\Sigma _{n}.$
\ So if we denote by $F^{\prime }$ the fixed set of $S$ in $U(n)/T^{n}$ we
see that $H_{S}^{\ast }(F^{\prime })$ is again, in a natural way, the
regular representation of $\Sigma _{n}$ over $\Lambda .$

The map $f:C_{n}(R^{3})\rightarrow U(n)/T^{n}$ of \cite{1}, being compatible
with $S$ must map $F$ to $F^{\prime }$ and since it is also compatible with $%
\Sigma _{n}$ it identifies the components of $F$ with the points of $%
F^{\prime }.$ \ The induced map $f^{\ast }$ on equivariant cohomology 
\begin{equation*}
H_{S}^{\ast }(F^{\prime })\rightarrow H_{S}^{\ast }(F)
\end{equation*}
is thus the standard identification of these two copies of the regular
representation of $\Sigma _{n}$ over $\Lambda .$

Finally, since $U(n)/T^{n}$ has no odd-dimensional cohomology, $H_{S}^{\ast
}(U(n)/T^{n})$ is a free $\Lambda $-module generated by a lift of a basis of 
$H^{\ast }(U(n)/T^{n}).$

Consider now the commutative diagram induced by $f.$%
\begin{equation*}
\begin{array}{ccc}
H_{S}^{\ast }(U(n)/T^{n}) & \overset{f^{\ast }}{\rightarrow } & H_{S}^{\ast
}(C_{n}(R^{3})) \\ 
\downarrow &  & \downarrow \\ 
H_{S}^{\ast }(F^{\prime }) & \widetilde{\rightarrow } & H_{S}^{\ast }(F)
\end{array}
.
\end{equation*}

\noindent By (2.1) the vertical arrows are isomorphisms modulo $t$-torsion.
\ Since all $\Lambda $ modules in this diagram are free it follows that the
vertical arrows and $f^{\ast }\mathbf{\ }$\textbf{are injective with a }$t$-%
\textbf{torsion cokernel.}

This proves the analogue of Theorem 1 with $SO(3)$ replaced by $S=SO(2).$ \
But, over $Q,$ the $SO(3)$-equivariant cohomology is just the subspace of $S$%
-equivariant cohomology invariant under $t\rightarrow -t$ (with $u=t^{2}).$
Then Theorem 1 follows from the case we have proved. \ Essentially this
means that the matrix $A(t)$ describing $f^{\ast }$ has no odd powers of $t,$
and so we can put $t^{2}=u.$

Note that Theorem 1, involving $SO(3),$ does not involve breaking the
symmetry. \ In the proof we picked an axis and worked with $SO(2),$ but this
was for convenience only. \ If we had worked with $SO(3)$ directly then
fixed points would have been replaced by $2$-dimensional orbits.

Now as promised we will discuss the question of lifting generators. \
Actually we are not interested in individual generators, but in a generating
subspace. \ Thus if $H_{S}^{\ast }(X)\rightarrow H^{\ast }(X)$ is surjective
we want to find a natural right inverse $H^{\ast }(X)\rightarrow H_{S}^{\ast
}(X).$ \ When a group like $\Sigma _{n}$ acts on both sides we want to
choose this right inverse to be compatible with the action of $\Sigma _{n}.$
\ Once this has been done we will get a natural isomorphism 
\begin{equation*}
H^{\ast }(X)\otimes \Lambda \rightarrow H_{S}^{\ast }(X)
\end{equation*}
compatible with $\Sigma _{n}.$

One way to pick such a right inverse is if there are invariant metrics. \ We
then pick the orthogonal complement to the kernel.

We begin by considering the case $X=U(n)/T^{n}$ with the full action of $%
U(n) $ (not just the subgroup $S)$ we then have the fibration 
\begin{equation*}
BT^{n}\rightarrow BU(n)
\end{equation*}
and the cohomology of the total space is the polynomial algebra in $n$
variables $t_{1},...,t_{n}.$ \ On this (in each degree) we can pick a $%
\Sigma _{n}$-invariant metric (e.g. the natural metric induced by a metric
on $C^{n})$. \ Taking the orthogonal complement then defines a natural right
inverse of the restriction 
\begin{equation*}
H^{\ast }(B_{T^{\ast }})\rightarrow H^{\ast }(U(n)/T^{n})
\end{equation*}
i.e. of the homomorphism 
\begin{equation*}
H_{U(n)}^{\ast }(U(n)/T^{n})\rightarrow H^{\ast }(U(n)/T^{n}).
\end{equation*}
This then gives a lift for all subgroups of $U(n)$ and so in particular for $%
SO(3)$ or its maximal torus $S.$ \ Since $\Sigma _{n}$ acts just by
permuting $t_{1},...,t_{n}$ it preserves the invariant metric and hence our
chosen lift.

Consider next the case $X=C_{n}(R^{3})$ and introduce the space 
\begin{equation*}
Y=\dprod\limits_{i,j}Y_{ij}
\end{equation*}
where each $Y_{ij}$ is a copy of the $S^{2}$ of directions in $R^{3}.$ \
Then we have maps 
\begin{equation*}
\alpha _{ij}:X\rightarrow Y_{ij}
\end{equation*}
by taking the direction of the line $(x_{i},x_{j})$ where $%
x=(x_{1},...,x_{n})\in X=C_{n}(R^{3}).$ \ Putting these together we get 
\begin{equation*}
\alpha :X\rightarrow Y
\end{equation*}
which is compatible with the natural actions of $\Sigma _{n}$ and $SO\left(
3\right) $ on both sides. \ For each $2$-sphere $Y_{ij}$ the tangent $SO(2)$%
-bundle $L$ is an $S$-bundle whose Chern class $c_{1}(L)$ is twice the
generator $y_{ij}\in H^{2}(Y,Z).$ \ Thus $\dfrac{1}{2}c_{1}(L)$ is an $S$%
-equivariant class which lifts $y_{ij}.$ \ This procedure gives us a natural
lift for each 
\begin{equation*}
y_{ij}\in H^{2}(Y).
\end{equation*}

Now we have a diagram 
\begin{equation*}
\begin{array}{ccc}
H_{S}^{\ast }(Y) & \overset{\dashleftarrow }{\rightarrow } & H^{\ast }(Y) \\ 
\downarrow &  & \downarrow \\ 
H_{S}^{\ast }(X) & \rightarrow & H^{\ast }(X)
\end{array}
\end{equation*}
in which the dotted line is the lift we have just described. \ To get a lift
for $X$ we note that $H^{\ast }(Y)\rightarrow H^{\ast }(X)$ is surjective 
\cite{1} so we only need to find a natural orthogonal to the kernel. \ This
will give a lift from $H^{\ast }(X)$ to $H^{\ast }(Y)$ and then we follow
with the lift to $H_{S}^{\ast }(Y)$ and project back to $H_{S}^{\ast }(X).$
\ Elementary diagram chasing shows that this gives what we want. \ To find
the orthogonal complement in $H^{\ast }(Y)$ we need a natural metric but
this is easy since 
\begin{equation*}
H^{\ast }(Y)=\bigotimes_{i,j}H^{\ast }(Y_{ij})
\end{equation*}
and each factor can be given a natural metric (since $Y_{ij}=S^{2}).$ \
Since no arbitrary choices were involved and since all indices $1,...,n$
were treated on an equal footing it follows that our lift is compatible with
the action of $\Sigma _{n}.$

Putting all this together we see finally that, as indicated in \S 1, Theorem
1 leads to a canonical homomorphism 
\begin{equation*}
A(u):H^{\ast }(U(n)/T^{n})\rightarrow H^{\ast }(C_{n}(R^{3}))\otimes Q[u]
\end{equation*}
where

\begin{equation*}
A(u)=A_{0}+A_{1}u+A_{2}u^{2}+...
\end{equation*}
is a polynomial in $u$ whose coefficients are homomorphisms 
\begin{equation*}
A_{k}:H^{\ast }(U(n)/T^{n})\rightarrow H^{\ast }(C_{n}(R^{3}))
\end{equation*}
which lower degree by $4k,$ and commute with the action of $\Sigma _{n}.$ \
For $k=0,$ $A_{0}$ is just the original homomorphism $f^{\ast },$ the
interest lies in the higher $A_{k}.$

\section{Grassmannians}

Theorem 1 can be generalised when we replace the flag manifold by a
Grassmannian and modify our configuration space accordingly. \ We proceed as
follows.

First we define a new configuration space $C_{r,s}(R^{3}).$ \ This
parametrizes two \textbf{unordered} sets of points of $R^{3}%
\;(x_{1},...,x_{r})$ and $(y_{1},...,y_{s})$ where $x_{i}\neq y_{j}$ for any 
$i,j.$ However we allow the $x_{i}$ to coincide among themselves and
similarly the $y_{j}$ . \ Clearly we have a natural map 
\begin{equation*}
C_{n}(R^{3})\rightarrow C_{r,s}(R^{3})\quad n=r+s
\end{equation*}

\noindent where $x_{1},...,x_{r}$ are the first $r$ points of a
configuration in $C_{n}(R^{3}).$

We now look carefully at the map 
\begin{equation*}
f_{n}:C_{n}(R^{3})\rightarrow U(n)/T^{n}
\end{equation*}
constructed in \cite{1} to see what happens when two of the points of the
configuration come together. \ Suppose $x_{2}\rightarrow x_{1}.$ \ The key
part of the construction in \cite{1} is that we define points $t_{ij}\in
S^{2}$ and then polynomials $p_{i}$ whose roots are the $t_{ij}$ (for $j\neq
i).$ \ As $x_{2}\rightarrow x_{1}$ the roots $t_{12}$ and $t_{21}$ cease to
be well-defined in the limit. \ However all other $t_{ij}$ are well-defined.
\ Thus the individual polynomials $p_{1}$ and $p_{2}$ are ill-defined but 
\textbf{the linear space spanned by them remains well-defined:} it consists
of all polynomials of degree $n-1$ which have $t_{1j}$ \ $(j\neq 1,2)$ as $%
n-2$ of their roots.

Similar remarks apply if several points coincide, so that by continuity we
can associate to any configuration $(x,y)\in C_{r,s}(R^{3})$ two linear
subspaces $C^{r}$ and $C^{s}$ of $C^{n}.$ \ Again by continuity (and the
properties of $f_{n})$ it follows that $C^{r}$ and $C^{s}\,$span$\,C^{n}.$ \
The orthogonalisation process used in \cite{1}, depending on the polar
decomposition, then gives us a map 
\begin{equation*}
f_{r,s}=C_{r,s}(R^{3})\rightarrow U(n)/U(r)\times U(s)\quad n=r+s
\end{equation*}
and a commutative diagram 
\begin{equation*}
\begin{array}{ccc}
C_{n}(R^{3}) & \overset{f_{n}}{\rightarrow } & U(n)/T^{n} \\ 
\downarrow &  & \downarrow \\ 
C_{r,s}(R^{3}) & \overset{f_{r,s}}{\rightarrow } & U(n)/U(r)\times U(s)
\end{array}
\end{equation*}
all of which is compatible with the action of $SO(3).$

We now consider the induced homomorphism $f_{r,s}^{\ast }$ in $SO(3)$%
-equivariant cohomology. \ First we restrict to $S=SO(2)\subset SO(3)$ and
consider the fixed point sets in the two spaces.

For $C_{r,s}(R^{3})$ the fixed point set consists of $C_{r,s}(R^{1}).$ \
Since we ignore the orderings of the $x_{i}$ and of the $y_{j}$ there are 
\begin{equation*}
\dfrac{n!}{r!s!}\text{ components}
\end{equation*}
Each component is contractible (since each $x_{i}$ can move in the interval
between the nearest $y_{j}$ and vice-versa).

For the Grassmannian the fixed points are finite in number and they are just
the images of the $n!$ fixed points in the flag manifold under projection
which acts by factoring out $\Sigma _{r}\times \Sigma _{s},$ giving 
\begin{equation*}
\frac{n!}{r!s!}\text{ fixed points.}
\end{equation*}

From the diagram it is easy to see that $f_{r,s}$ sends each component of $%
C_{r,s}(R^{1})$ to a different fixed point. \ Exactly as in the proof of
Theorem 1 we then get.\medskip

\noindent \textbf{Theorem 2.\quad The map }$f_{r,s}:C_{r,s}(R^{3})%
\rightarrow U(n)/U(r)\times U(s))$ \textbf{induces a homomorphism} 
\begin{equation*}
f_{r,s}^{\ast }:H_{SO(3)}^{\ast }(U(n)/U(r)\times U(s),Q)\rightarrow
H_{SO(3)}^{\ast }(C_{r,s}(R^{3}),Q)
\end{equation*}
\textbf{which is injective and has cokernel annihilated by a power of }$%
u.\medskip $

It is clear that Theorem 2 can be further generalised to any other
homogeneous space of $U(n)$ which occurs as an adjoint orbit. \ These are
all quotients by subgroups of the form 
\begin{equation*}
U(r_{1})\times U(r_{2})\times ...\times U(r_{t})\quad \Sigma r_{i}=n.
\end{equation*}
Moreover if any of the $r_{i}$ are equal there will be a small symmetric
group permuting these blocks and the corresponding map from the generalised
configuration space will commute with this. \ Theorem 1 is the special case
where all $r_{i}=1.$

Note that if $r=1,\;s=n-1$ then $C_{1,n-1}(R^{3})$ has the homotopy type of
the $(n-1)^{\text{th}}$ symmetric product of $S^{2},$ i.e. the complex
projective space $P_{n-1}$ and $f_{1,n-1}^{\ast }$ in Theorem 2 is then
obviously an isomorphism

\section{K-Theory}

We shall essentially repeat the arguments of \S 2 with $K$-theory replacing
cohomology. \ Thus we want to consider the homomorphism in $K_{S}^{\ast }$
induced by our map $f$: 
\begin{equation*}
f^{\ast }:K_{S}^{\ast }(U(n)/T^{n})\rightarrow K_{S}^{\ast }(C_{n}(R^{3}))
\end{equation*}
where as before $S=SO(2)$ is given by rotation about a preferred axis in $%
R^{3}$ (and $S$ acts on $U(n)/T^{n}$ through the principal $SU(2)$ subgroup).

First we want to show that for both spaces $K_{S}^{1}=0$ and that $K_{S}^{0}$
is a free $R[S]$-module of $\func{rank}n!$ \ Note that our ground ring is
now the representation ring of the circle 
\begin{equation*}
R[S]=Z[q,q^{-1}].
\end{equation*}
For brevity we denote this by $A.$

Consider first the flag manifold 
\begin{equation*}
X=U(n)/T^{n}
\end{equation*}
and recall that, as a complex manifold (depending on an ordering of $%
1,...,n),$ it has a decomposition into even dimensional cells (indexed by $%
\Sigma _{n}).$ \ These (Bruhat) cells are the orbits of a Borel subgroup $%
B\subset GL(n,C)$ (the upper triangular matrices). \ In particular they are
preserved by the maximal torus $T^{n}$ of $U(n)$ and hence by the circle $%
S\subset T^{n}.$ \ Moreover each cell contains just one fixed point of $S.$
\ The union of cells of dimension $\leq 2p$ gives a closed subspace $%
X_{p}\subset X$ and a standard induction argument using the $X_{p}$ enables
one to calculate the cohomology or the $K$-theory of $X.$ \ Since $S$
preserves the $X_{p}$ the same arguments yield the calculation of $%
K_{S}^{\ast }(X)$ showing that $K_{S}^{1}(X)=0$ and that $K_{S}^{0}(X)$ is a
free $A$-module of $\func{rank}n!$ \ 

Using some algebraic geometry one can even get a natural basis indexed by
the cells, i.e. by $\Sigma _{n}.$ \ The closure of each cell is an algebraic
subvariety and its sheaf (via a resolution) determines an element of the $K$%
-group (the algebraic and topological $K$-groups of the flag manifold
coincide). \ Since everything is invariant under $T^{n}$ and hence under $S$
we get a basis of $K_{S}$ of the flag manifold. \ Note however that this
basis is definitely not compatible with the action $\Sigma _{n}$ (since it
depended on an ordering).

For $X=C_{n}(R^{3})$ we have no such convenient cell decomposition so we
proceed differently. \ Recall that $S$ has $n!$ contractible fixed
components (indexed by $\Sigma _{n})$ and all other orbits are actually
free. \ It follows from general localisation theory for equivariant $K$%
-theory \cite{7} that $K_{S}^{\ast }(X)$ can only have $A$-torsion at the
identity, i.e. annihilated by a power of $(q-1).$ \ On the other hand a
general theorem on completion of equivariant $K$-theory\ \cite{3} asserts
that 
\begin{equation*}
K_{G}^{\ast }(X)^{\wedge }\cong K(X_{G}^{\ast })
\end{equation*}
where $\wedge $ denotes completion at the identity, i.e. in the $I(G)$-adic
topology. \ For $G=S$ this means the $(q-1)$-adic topology. \ Since this
does not kill torsion at the identity any torsion in $K_{S}^{\ast }(X)$
would persist to $K^{\ast }(X_{S}),$ but $X_{S}$ is the fibration over $%
CP_{\infty }$ with fibre $X$ and, just as in \S 2, a spectral sequence
argument (together with the vanishing of all odd cohomology of $X)$ then
shows that $K^{1}(X_{S})=0$ and $K^{\ast }(X_{S})$ is a free $\hat{A}$%
-module. \ This shows that there could not have been any $A$-torsion in $%
K_{S}^{0}(X)$ and that $K_{S}^{1}(X)=0.$

Note that this argument would not have worked if there were $1$-dimensional
orbits with non-trivial isotropy group. \ We might then have had torsion
supported at points other than $q=1$ and this would not have been detected
in the completion. \ Such orbits actually occur for the action on the flag
manifold even though, as our other argument shows, there is no torsion. \
However such orbits will affect the torsion statement in Theorem 3 as we
shall see.

Just as in \S 2 we now consider the restriction to the fixed point sets of $%
S $ in the two spaces. \ The localisation theorems for $K_{S}^{\ast }$ \cite
{7} tell us that the restrictions are isomorphisms modulo $A$-torsion. \
Hence, as before, we conclude

\noindent \textbf{Theorem 3.\quad The homomorphism} 
\begin{equation*}
f^{\ast }:K_{S}^{\ast }(U(n)/T^{n})\rightarrow K_{S}^{\ast }(C_{n}(R^{3}))
\end{equation*}
\textbf{is injective and its cokernel is a torsion module over the ring }$%
R[S]=Z[q,q^{-1}].$ \ \textbf{Here} $S=SO(2).\bigskip $

\noindent \textbf{Remarks}

\begin{enumerate}
\item  We can be more precise about the torsion module. \ The localisation
theorem tells us that its support is contained in the set of all $q$ which
are $m^{\text{th}}$ roots of unity, where $m$ runs over the orders of the
finite isotropy groups of $S$ acting on $U(n)/T^{n}.$ \ These are just the
values $m\leq n-1.$

\item  With a bit more work it is also possible to formulate Theorem 3 with $%
S$ replaced by $SO(3)$ or $SO(2).$

\item  Theorem 3 also generalises to other Lie groups as will be shown on a
later occasion.$^{{}}$

\item  The homomorphism $f^{\ast }$ of Theorem 3 is of course a ring
homomorphism not just an $A$-module homomorphism. \ The ring structures of
both sides can be explicitly described and $f^{\ast }$ can then be described
in these terms. \ Similar remarks apply to Theorems 1 and 2 and I\ hope to
follow these up elsewhere. \ My thanks are due to William Graham for
interesting me in this aspect.\vspace{0.5in}
\end{enumerate}

\vspace{1in}

\noindent Department of Mathematics \& Statistics,

\noindent University of Edinburgh,

\noindent James Clerk Maxwell Buildings,

\noindent King's Buildings,

\noindent Edinburgh \ EH9 3JZ.

\end{document}